\documentclass{article}
\usepackage{latexsym,amssymb,graphicx}
\input amssym.def \input amssym

\newtheorem{te}{Theorem}[section]

\newtheorem{de}[te]{Definition}
\newtheorem{lm}[te]{Lemma}

\newtheorem{pp}[te]{Proposition}

\newtheorem{ex}[te]{Example}

\def\dokaz{\noindent{\bf Proof. }}
\def\kraj{\hfill $\Box$ \par \vspace*{2mm} }

\def\widemid{\hspace{1mm}\widetilde{\mid}\hspace{1mm}}

\newcommand{\zve}[1]{{{}^*\hspace{-0.5mm}#1}}
\newcommand{\zvez}[1]{{{}^*\hspace{-1mm}#1}}
\def\zvezmid{\hspace{1mm}\zvez\mid\hspace{1mm}}
\def\zs{\forall}
\def\po{\exists}
\def\str{\rightarrow}
\def\Str{\Rightarrow}

\def\dl{\Leftrightarrow}
\def\gstr{\hspace{-0.1cm}\uparrow\hspace{0.1cm}}
\def\dstr{\hspace{-0.1cm}\downarrow}
\def\ps{\subseteq}

\def\cU{{\cal U}}
\def\cV{{\cal V}}
\def\cW{{\cal W}}

\def\cF{{\cal F}}
\def\cG{{\cal G}}
\def\cH{{\cal H}}
\def\cP{{\cal P}}
\def\cQ{{\cal Q}}
\def\cM{{\cal M}}

\begin{document}
\begin{center}
           {\huge \bf More number theory in $\beta N$}
\end{center}
\begin{center}
{\small \bf Boris  \v Sobot}\\[2mm]
{\small  Faculty of Sciences, University of Novi Sad,\\
Trg Dositeja Obradovi\'ca 4, 21000 Novi Sad, Serbia\\
e-mail: sobot@dmi.uns.ac.rs, ORCID: 0000-0002-4848-0678}
\end{center}
\begin{abstract} \noindent
We continue the research of an extension $\widemid$ of the divisibility relation to the Stone-\v Cech compactification $\beta N$. First we prove that ultrafilters we call prime actually possess the algebraic property of primality. Several questions concerning the connection between divisibilities in $\beta N$ and nonstandard extensions of $N$ are answered, providing a few more equivalent conditions for divisibility in $\beta N$. Results on uncountable chains in $(\beta N,\widemid)$ are proved and used in a construction of a well-ordered chain of maximal cardinality. Finally, we consider ultrafilters without divisors in $N$ and among them find the maximal class.

\vspace{1mm}


\vspace{1mm}

\end{abstract}

\section{Introduction}

For any set $S$ let $\beta S$ denote the set of ultrafilters on $S$. For each $n\in S$ the principal ultrafilter $\{A\subseteq S:n\in A\}$ is identified with $n$, so $\beta S$ is thought of as containing $S$. The ultrafilters in $\beta S\setminus S$ are called {\it free} (nonprincipal). If $S$ is endowed with discrete  topology, $\beta S$ (with base sets $\overline{A}=\{\cF\in\beta S:A\in\cF\}$) is known as the Stone-\v Cech compactification of this space.

Every semigroup $(S,*)$ with discrete topology can be extended to a right-topological semigroup on $\beta S$. This extension has many nice properties which can be used to obtain various interesting results about $(S,*)$. Many such results can be found in \cite{HS}.

One of the most natural examples is $(N,\cdot)$ - the set of natural numbers with multiplication. Here the operation on $\beta N$ is defined with
$$A\in\cF\cdot\cG\dl\{n\in S:A/n\in\cG\}\in\cF,$$
where $A/n=\{m\in N:mn\in A\}$. Having extended multiplication, one can ask: how to extend the divisibility relation? There seem to be more than one natural way to define the extension; four possible divisibility relations on $\beta N$ were first examined in \cite{So1}. Each of these relations has the usual $\mid$ as the restriction to $N\times N$, and moreover they all coincide on $N\times\beta N$: $\cF\in\beta N$ is divisible by $n\in N$ if and only if $nN:=\{m\in N:n\mid m\}\in\cF$.

The divisibility relation on $\beta N$ that has the nicest properties is the following: for $\cF,\cG\in\beta N$,
$$\cF\widemid\cG\dl\cF\cap\cU\subseteq\cG,$$
where $\cU=\{A\subseteq N:A\gstr=A\}$ and $A\gstr=\{n\in N:\po a\in A\;a\mid n\}$. $\widemid$ is a quasiorder (reflexive and transitive). We will actually think of it (without explicitly mentioning it) as of the order obtained from $\widemid$ on $\beta N/=_\sim$, where $\cF=_\sim\cG$ if and only if $\cF\widemid\cG$ and $\cG\widemid\cF$. It is weaker than the other three relations; in particular, $\cG=\cH_1\cdot\cF\cdot\cH_2$ implies $\cF\widemid\cG$. Throughout the paper, whenever we speak of divisibility, it is understood that it is $\widemid$-divisibility. $[\cF]$ denotes the $=_\sim$-equivalence class of $\cF$.

The investigation of the structure of divisibility hierarchy of $\beta N$ began with the idea to eventually apply it to $(N,\mid)$. So far we have seen that many properties of the divisibility on $N$ reflect in some way to $\beta N$; we recapitulate some of them here. Our hope is that we will be able to acquire results in the other direction, as it was done with extensions of operations.

In the $\widemid$-hierarchy, the ultrafilters of $\beta N$ (or, rather, their equivalence classes modulo $=_\sim$) are divided into two parts. The first, "lower" part, consists of $\omega$-many levels $\overline{L_n}$, resembling in many ways the order $(N,\mid)$, and including it. All the $=_\sim$-equivalence classes here are sigletons (\cite{So3}, Lemma 5.13). So, $1$ is the smallest element. The level $\overline{L_1}$ consists of ultrafilters containing the set $P=L_1$ of prime numbers; we also call these ultrafilters {\it prime}. There are $2^{\goth c}$ prime ultrafilters. Section 2 of this paper gives us one more of their interesting properties.

In general, the $n$-th level $\overline{L_n}$ consists of ultrafilters containing the set $L_n$ of numbers having exactly $n$ (not necessarily distinct) prime factors. There is no prime factorization theorem here, but we are able to decompose every ultrafilter from $L=\bigcup_{n<\omega}\overline{L_n}$ into basic ingredients. To make the description more precise, we use as basic not only prime ultrafilters but also their powers: $\cP^n$ is generated by $\{A^n:A\in\cP,A\subseteq P\}$, where $A^n=\{p^n:p\in A\}$. So each ultrafilter from $\overline{L_n}$ can be decomposed into exactly $n$ parts, with $\cP^n$ counted $n$ times, and same basic parts are allowed to occur multiple times. 

For example, there are three types of ultrafilters on the second level: the squares $\cP^2$ (exactly one for each prime $\cP$), ultrafilters having two distinct primes below them, but also ultrafilters "divisible twice" by a prime. More precisely, there is a prime ultrafilter having $2^{\goth c}$ successors in $\overline{L_2}$ (\cite{So3}, Theorem 3.13). To understand this phenomenon better, we needed to make a connection of $\beta N$ with a nonstandard universe $V(\zve N)$. (We explain some basic nonstandard notions in a separate subsection below.) Each $\cF\in\beta N$ has a corresponding set of nonstandard integers $\mu(\cF)$ in $\zve N$, called the monad of $\cF$. If the extension is an enlargement, all monads are nonempty. In \cite{So4}, Theorem 3.1, we established an important connection between divisibility relations $\zvezmid$ in $V(\zve N)$ and $\widemid$ in $\beta N$: for $\cF,\cG\in\beta N$, $\cF\widemid\cG$ holds if and only if there are $x\in\mu(\cF)$ and $y\in\mu(\cG)$ such that $x\zvezmid y$. Moreover, a nonstandard integer is on the $n$-th level of the $\zvezmid$-hierarchy if and only if it belongs to the monad of an ultrafilter from $\overline{L_n}$. Squares of primes in $\zve N$ are in the monads of squares of prime ultrafilters; products of two nonstandard primes belonging to distinct monads generate ultrafilters $\widemid$-divisible by two distinct prime ultrafilters; finally, products of nonstandard primes belonging to the same monad generate ultrafilters $\widemid$-divisible by only one prime ultrafilter.

More results about the "lower" part of the $\widemid$-hierarchy can be found in \cite{So3}. What about the "upper" part? Things get more complicated there, and ultrafilters can no longer be organized by levels (see \cite{So4}). It seems that representing these ultrafilters as limits of $\widemid$-ascending chains can be the right way to consider them. We recapitulate what we know about these limits in another separate subsection.

Finally, on the very top of the hierarchy there is the greatest class MAX, consisting of ultrafilters divisible by all others. Their existence is easy to show, since $\cU$ has the finite intersection property. In a subsequent paper we will try to find the place of MAX among some other classes of ultrafilters, important for topological dynamics. In the last section of this paper we will encounter one more interesting $=_\sim$-class, NMAX.\\

{\bf Notation.} $N$ is the set of natural numbers (without zero) and $P$ denotes the set of (standard) prime numbers. To make our statements easier to read, we reserve calligraphic letters $\cF,\cG,\cH,\dots$ for ultrafilters, with $\cP,\cQ,\dots$ denoting prime ultrafilters. Small letters $x,y,z,\dots$ will denote elements of $\zve N$, with $p,q,\dots$ reserved for primes. For $A,B\subseteq N$, $A^c=N\setminus A$, $A^2=\{a^2:a\in A\}$ (to avoid confusion, we will not abbreviate $A\times A$ to $A^2$), $A^{(2)}=\{a_1a_2:a_1,a_2\in A\land a_1\neq a_2\}$ and $AB=\{ab:a\in A\land b\in B\land a\neq b\}$. If also $n\in N$, then $nN=\{nm:m\in N\}$ and $A/n=\{\frac an:a\in A\land n\mid a\}$. We recall that $\cU=\{A\subseteq N:A\gstr=A\}$, where $A\gstr=\{n\in N:\po a\in A\;a\mid n\}$ and $\cV=\{A\subseteq N:A\dstr=A\}$, where $A\dstr=\{n\in N:\po a\in A\;n\mid a\}$. In particular, we write $a\gstr$ instead of $\{a\}\gstr$.\\

{\bf Nonstandard methods.} We follow the superstructure approach of Robinson and Zakon from \cite{RZ}. Let $X$ be a set containing (a copy of) $N$. We assume that elements of $X$ are atoms: none of them contains as an element any of the others. Let $V_0(X)=X$, $V_{n+1}(X)=V_n(X)\cup P(V_n(X))$ for $n\in\omega$ and $V(X)=\bigcup_{n<\omega}V_n(X)$. $V(X)$ is called a {\it superstructure}. Superstructures are convenient because they include mostly everything one needs when working with $X$: subsets, relations, functions... For example, the divisibility relation $\mid$ is in $V_3(N)$. If $V(X)$ is a superstructure, its {\it nonstandard extension} is a pair $(V(Y),*)$, where $V(Y)$ is a superstructure with the set of atoms $Y\supset X$ and $*:V(X)\rightarrow V(Y)$ is a rank-preserving function such that $\zve X=Y$ and satisfying the following principle.

{\it The Transfer Principle.} For every bounded formula $\varphi$ and every $a_1,a_2,\dots$, $a_n\in V(X)$, $\varphi(a_1,a_2,\dots,a_n)$ holds in $V(X)$ if and only if $\varphi(\zve a_1,\zve a_2,\dots,\zve a_n)$ holds in $V(Y)$.

(A first-order formula is bounded if all its quantifiers are bounded, i.e.\ of the form $(\forall x\in y)$ or $(\exists x\in y)$. The free variables that appear in $\varphi(a_1,a_2,\dots,a_n)$ are exactly objects $a_1,a_2,\dots,a_n$ from $V(X)$ and in $\varphi(\zve a_1,\zve a_2,\dots,\zve a_n)$ they are replaced with their star-counterparts. The atomic subformulas in $\varphi$ are of the form $A(x_1,\dots,x_k)$ for some $k$-ary relation $A\in V(X)$, which also gets replaced with $\zve A$. For example, if the divisibility relation $\mid$ appears in $\varphi$, it gets replaced with $\zve\mid$. Exceptions to this rule can be introduced, so one doesn't write $\zve =$ instead of $=$; see \cite{H} for justification.)

Since all objects we need live in the superstructure $V(N)$, we will consider only nonstandard extensions of this superstructure. Often we call just the set $V(\zve N)$ the nonstandard extension of $V(N)$. A set in $V(\zve N)$ is {\it internal} if it belongs to $\zve A$ for some $A\in V(N)$. So quantifiers in formulas applicable for the Transfer Principle range only through internal sets. {\it The Internal Definition Principle} claims that, basically, sets defined from internal sets are also internal.

We will often need our nonstandard extensions to satisfy additional conditions. We call $V(\zve N)$ a $\kappa$-{\it enlargement} if for every family $F$ of subsets of some set in $V(N)$ with the finite intersection property such that $|F|<\kappa$ there is an element in $\bigcap_{A\in F}\zve A$. (This definition, more common in recent papers, is somewhat different from the one we used in \cite{So4}, but the proofs thereof easily translate.) $V(\zve N)$ is $\kappa$-{\it saturated} if every family $F$ of internal sets in $V(\zve N)$ with the finite intersection property such that $|F|<\kappa$ has nonempty intersection. $\kappa$-saturated extensions (and therefore also $\kappa$-enlargements) are known to exist in ZFC.

Monads were first introduced in \cite{Lu} and they present the connection between nonstandard universe and the Stone-\v Cech compactification. The {\it monad} of an ultrafilter $\cF\in\beta N$ is the set $\{x\in\zve N:(\zs A\in\cF)x\in\zve A\}$. For $x\in\zve N$, $v(x)$ is the unique ultrafilter $\cF\in\beta N$ such that $x\in\mu(\cF)$ (we call $v(x)$ the ultrafilter generated by $x$). If $V(\zve N)$ is a ${\goth c}^+$-enlargement, all monads of free ultrafilters are nonempty (and actually of the same cardinality as $\zve N$), giving us rich enough structure to work with. The stronger condition of ${\goth c}^+$-saturation provides us with more additional properties that we will investigate in this paper. For more information on nonstandard extensions the reader can consult \cite{CK}, section 4.5, \cite{H} or \cite{L}. In particular, \cite{L} contains some results, called the bridge theorems, that we will use in Section 3 to answer some questions from \cite{So4} and extend the list of equivalent conditions for $\widemid$-divisibility.\\

{\bf Chains.} In $\beta N$ only an eventually constant sequence can have a (standard) topological limit, so it is common to work with limits by ultrafilters. If $\langle\cG_i:i\in I\rangle$ is a sequence of ultrafilters on $N$ and $\cF$ is an ultrafilter on $I$, $\cG:=\lim_{i\str\cF}\cG_i$ is the ultrafilter on $N$ such that $A\in\cG$ if and only if $\{i\in I:A\in\cG_i\}\in\cF$. See also section 3.5 of \cite{HS}.

In \cite{So4} we considered $\widemid$-chains of ultrafilters of length $\omega$, and showed that every such chain has the smallest upper bound. In Section 4 of this paper we get more results about such (and longer) chains, hopefully getting closer to understanding the structure of "upper" half of the $\widemid$-hierarchy.

\section{More on prime ultrafilters}

In \cite{So1} we defined, for $\cF\in\beta N$, $D(\cF)=\{A\ps N:\{n\in N:A/n=N\}\in\cF\}=\{A\ps N:\{n\in N:nN\ps A\}\in\cF\}$. $D(\cF)$ is a filter contained in $\cF$ and $\cF\widemid\cG$ is equivalent to $D(\cF)\subseteq\cG$ (Theorem 6.2 of \cite{So1}). Note that
\begin{equation}
\cF\cap\cU\ps D(\cF)\label{eqDp}
\end{equation}
for every $\cF\in\beta N$. Namely, if $A\in\cF\cap\cU$ then for each $n\in A$ holds $nN\ps A$, i.e.\ $A\ps\{n\in N:nN\ps A\}$ so $\{n\in N:nN\ps A\}\in\cF$, which means that $A\in D(\cF)$. The reverse of (\ref{eqDp}) does not hold, since $D(\cF)\not\subseteq\cU$: if $A\in D(\cF)$ is arbitrary and $a\in A$ is not $\mid$-minimal, then $A\setminus\{a\}\notin\cU$, but $\{n\in N:nN\ps A\setminus\{a\}\}$ differs from $\{n\in N:nN\ps A\}$ by only finitely many elements, so it is also in $\cF$. Hence $A\setminus\{a\}\in D(\cF)\setminus\cU$.

Let $\cU'=\{B\gstr:B\ps P\}$. Clearly, $\cU'\ps\cU$.

\begin{lm}\label{prostilema}
Let $\cF,\cG\in\beta N$. Then:

(a) $D(\cF\cdot\cG)\cap\cU'\ps D(\cF)\cup D(\cG)$;

(b) $\cF\cdot\cG\cap\cU'\ps\cF\cup\cG$.
\end{lm}

\dokaz (a) Let $A=B\gstr$ for some $B\ps P$ and $A\in D(\cF\cdot\cG)$. We have:
\begin{eqnarray}
A\in D(\cF\cdot\cG) &\dl& \{n\in N:nN\ps A\}\in\cF\cdot\cG\nonumber\\
 &\dl& \{m\in N:\{n\in N:nN\ps A\}/m\in\cG\}\in\cF\nonumber\\
 &\dl& \{m\in N:\{n\in N:mnN\ps A\}\in\cG\}\in\cF\nonumber\\
 &\dl& \{m\in N:\{n\in N:nN\ps A/m\}\in\cG\}\in\cF\nonumber\\
 &\dl& \{m\in N:A/m\in D(\cG)\}\in\cF.\label{eqDy}
\end{eqnarray}
But $A/m=\left\{\begin{array}{ll}
N, & \mbox{if }p\mid m\mbox{ for some }p\in B\\
A, & \mbox{otherwise.}
\end{array}\right.$

So either $A\in D(\cG)$, or $\{m\in N:A/m\in D(\cG)\}=\{m\in N:A/m=N\}$, which means by (\ref{eqDy}) that $A\in D(\cF)$.\\

(b) Using (\ref{eqDp}) and (a) we have $\cF\cdot\cG\cap\cU'=(\cF\cdot\cG\cap\cU)\cap\cU'\ps D(\cF\cdot\cG)\cap\cU'\ps D(\cF)\cup D(\cG)\ps\cF\cup\cG$.\kraj

In the next theorem we prove that all ultrafilters in $\overline{P}$ deserve the name "prime" (in algebraic sense).

\begin{te}\label{prosti}
If $\cP\in\beta N$ is prime, then $\cP\widemid \cF\cdot\cG$ implies $\cP\widemid\cF$ or $\cP\widemid\cG$.
\end{te}

\dokaz Assume the opposite, that there are $A_\cF\in\cP\cap\cU\setminus\cF$ and $A_\cG\in\cP\cap\cU\setminus\cG$. Since $P\in\cP$, the sets $B_\cF=A_\cF\cap P$ and $B_\cG=A_\cG\cap P$ are also in $\cP$. $B_\cF\cap B_\cG\in\cP$ implies that the set $A=(B_\cF\cap B_\cG)\gstr$ is in $\cP\cap\cU'$, so by $\cP\widemid\cF\cdot\cG$ it is in $\cF\cdot\cG\cap\cU'$ as well. By Lemma \ref{prostilema}(b) it must also belong to $\cF$ or $\cG$, say $A\in\cF$. But $A\ps A_\cF$ and $A_\cF\notin\cF$, a contradiction.\kraj

\section{More about monads}

In \cite{L} several theorems were proved that enable us to translate formulas from $V(\zve N)$ to equivalent formulas in $V(N)$ and vice versa. The basic such theorem was called The Bridge Theorem there, so we will address all such results as bridge theorems. They can be thought of as a more comfortable way of applying the enlargement or saturation condition. In each of them $\phi(x_1,\dots,x_n)$ is a first-order formula and it is understood that $x_1,\dots,x_n$ are all free variables appearing in $\phi$. The first one we will need is Corollary 2.2.14 of \cite{L}.

\begin{pp}\label{bridge1}
Let $V(\zve N)$ be ${\goth c}^+$-saturated, $\cF,\cG\in\beta N$ and $z_1,\dots,z_k\in V(N)$. The following conditions are equivalent:

(i) $(\zs x_1,\dots,x_n\in\mu(\cF))(\po y_1,\dots,y_m\in\mu(\cG))\phi(x_1,\dots,x_n,y_1,\dots,y_m,\zve z_1$, $\dots,\zve z_k)$;

(ii) $(\zs B\in\cG)(\po A\in\cF)(\zs a_1,\dots,a_n\in A)(\po b_1,\dots,b_m\in B)\phi(a_1,\dots,a_n,b_1$, $\dots,b_m,z_1,\dots,z_k)$.
\end{pp}

We use it to answer (for ${\goth c}^+$-saturated extensions) affirmatively Question 5.3 left unresolved in \cite{So4}.

\begin{lm}\label{zspo}
Let $V(\zve N)$ be ${\goth c}^+$-saturated and $\cF,\cG\in\beta N$ are such that $\cF\widemid\cG$. Then:

(a) $(\zs x\in\mu(\cF))(\po y\in\mu(\cG))x\zvezmid y$;

(b) $(\zs y\in\mu(\cG))(\po x\in\mu(\cF))x\zvezmid y$.
\end{lm}

\dokaz (a) By Proposition \ref{bridge1} $(\zs x\in\mu(\cF))(\po y\in\mu(\cG))x\zvezmid y$ is equivalent to 
\begin{equation}\label{eqzspo}
(\zs B\in\cG)(\po A\in\cF)(\zs a\in A)(\po b\in B)a\mid b.
\end{equation}
So let $B\in\cG$. If $B\in\cF$, we let $A:=B$ and (\ref{eqzspo}) obviously holds. If $B\notin\cF$ let $B':=\{b\in B^c:b\gstr\subseteq B^c\}$ and $A:=B^c\setminus B'$. Then $B'\notin\cF$ (otherwise $B'\in\cF\cap\cU$ would imply $B'\in\cG$ but $B\cap B'=\emptyset$, a contradiction). So $A\in\cF$ and clearly (\ref{eqzspo}) holds again.\\

(b) is proven analogously, using $\cV$ in place of $\cU$.\kraj

Our next bridge theorem is similar to Theorem 2.2.9 from \cite{L}. For completeness' sake we include the proof here.

\begin{lm}\label{bridge2}
Let $V(\zve N)$ be a ${\goth c}^+$-enlargement, $\cF,\cG\in\beta N$ and $z_1,\dots,z_k\in V(N)$. The following conditions are equivalent:

(i) $(\po x_1,\dots,x_n\in\mu(\cF))(\po y_1,\dots,y_n\in\mu(\cG))\phi(x_1,\dots,x_n,y_1,\dots,y_m,\zve z_1,\dots$, $\zve z_k)$;

(ii) $(\zs A\in\cF)(\zs B\in\cG)(\po a_1,\dots,a_n\in A)(\po b_1,\dots,b_m\in B)\phi(a_1,\dots,a_n,b_1$, $\dots,b_m,z_1,\dots,z_k)$.
\end{lm}

\dokaz (i)$\Str$(ii) For each $A\in\cF$, $B\in\cG$, (i) implies $(\po x_1,\dots,x_n\in\zve A)(\po y_1,\dots$, $y_n\in\zve B)\phi(x_1,\dots,x_n,y_1,\dots,y_m,\zve z_1,\dots$, $\zve z_k)$ so, by transfer, $(\po a_1,\dots,a_n\in A)(\po b_1,\dots,b_m\in B)\phi(a_1,\dots,a_n,b_1,\dots,b_m,z_1,\dots,z_k)$.\\

(ii)$\Str$(i) Let
$$\Phi:=\{(a_1,\dots,a_n,b_1,\dots,b_m)\in N^{n+m}:\phi(a_1,\dots,a_n,b_1,\dots,b_m,z_1,\dots,z_k)\}.$$
Then
$$\zve\Phi=\{(x_1,\dots,x_n,y_1,\dots,y_m)\in(\zve N)^{n+m}:\phi(x_1,\dots,x_n,y_1,\dots,y_m,\zve z_1,\dots,\zve z_k)\}.$$
We prove that the family $\{\Phi\}\cup\{A^n\times B^m:A\in\cF,B\in\cG\}$ has the finite intersection property. Since $\cF$ and $\cG$ are closed for intersections, it is enough to see that each $A^n\times B^m$ intersects $\Phi$, which follows from (ii). By the ${\goth c}^+$-enlarging property, there is $(x_1,\dots,x_n,y_1,\dots$, $y_m)\in\zve\Phi\cap\bigcap_{A\in\cF,B\in\cG}(\zve A)^n\times(\zve B)^m$. Then $x_i\in\mu(\cF)$ and $y_j\in\mu(\cG)$ for all $i,j$. So (i) holds.\kraj

Now we can add several more equivalent conditions to Theorem 3.1 from \cite{So4} and obtain a better view of the connection between divisibilities in $V(\zve N)$ and $\beta N$.

\begin{te}\label{ekviv}
The following conditions are equivalent for every two ultrafilters $\cF,\cG\in\beta N$:

(i) $\cF\widemid\cG$;

(ii) $(\zs A\in\cF)(\zs B\in\cG)(\po a\in A)(\po b\in B)a\mid b$;

(iii) in every ${\goth c}^+$-enlargement $V(\zve N)$, there are $x\in\mu(\cF)$, $y\in\mu(\cG)$ such that $x\zvezmid y$;

(iv) in some ${\goth c}^+$-enlargement $V(\zve N)$, there are $x\in\mu(\cF)$, $y\in\mu(\cG)$ such that $x\zvezmid y$;

(v) in every ${\goth c}^+$-saturated extension $V(\zve N)$, for every $x\in\mu(\cF)$ there is $y\in\mu(\cG)$ such that $x\zvezmid y$;

(vi) in every ${\goth c}^+$-saturated extension $V(\zve N)$, for every $y\in\mu(\cG)$ there is $x\in\mu(\cF)$ such that $x\zvezmid y$.
\end{te}

\dokaz (i)$\Str$(ii) If $A\in\cF$ and $B\in\cG$, (i) implies that $A\gstr\in\cG$, so $B\cap A\gstr\in\cG$ as well. If $b\in B\cap A\gstr$ then there is $a\in A$ such that $a\mid b$ and we are done.\\


(ii)$\Str$(iii) follows directly from Lemma \ref{bridge2}.\\

(i)$\dl$(iii)$\dl$(iv) was proved in \cite{So4}, Theorem 3.1.\\

(i)$\Str$(v)$\land$(vi) follows from Lemma \ref{zspo}.\\

(v)$\Str$(iv) and (vi)$\Str$(iv) are trivial, since every ${\goth c}^+$-saturated extension is a ${\goth c}^+$-enlargement.\kraj

We note that extensions of relations from a set $X$ to $\beta X$ were considered in general in literature. The reason for including (ii) as a separate condition in the theorem above is that it shows that our relation $\widemid$ is exactly what is called "canonical extension" in \cite{PS}. This is another argument showing that $\widemid$ may be "the right" divisibility relation to consider on $\beta N$.\\

Let us call a set $X\subseteq\zve N$ {\it convex} if for all $x,y\in X$ and $z\in\zve N$, $x\zvezmid z$ and $z\zvezmid y$ imply $z\in X$. The following result answers negatively Question 5.4 from \cite{So4}.

\begin{lm}
Let $V(\zve N)$ be a ${\goth c}^+$-enlargement. For every $\cF\in MAX$, $\mu(\cF)$ is not a convex set.
\end{lm}

\dokaz Let $\cG\neq\cF$ be any other ultrafilter in $MAX$. We want to prove that $(\po x,y\in\mu(\cF))(\po z\in\mu(\cG))(x\neq y\neq z\land x\zvezmid z\land z\zvezmid y)$. By Lemma \ref{bridge2} this is equivalent to
\begin{equation}\label{eqpopo}
(\zs A\in\cF)(\zs C\in\cG)(\po a,b\in A)(\po c\in C)(a\neq c\neq b\land a\mid c\land c\mid b).
\end{equation}
Let $A\in\cF$ and $C\in\cG$, and let $a\in A$ be arbitrary. Since $aN\in\cG$, $aN\cap C$ is infinite so there is $c\in C\setminus\{a\}$ such that $a\mid c$. Analogously, since $cN\in\cF$, $cN\cap A$ is infinite so there is $b\in A\setminus\{c\}$ such that $c\mid b$. Hence (\ref{eqpopo}) holds.\kraj

There are, on the other hand, ultrafilters $\cF$ such that $\mu(\cF)$ is convex: if $\cF$ contains an infinite antichain $A$ (for example, if $P\in\cF$) then obviously (\ref{eqpopo}) does not hold for that choice of $A$.

The following (embarrassingly simple) example answers negatively Questions 5.1 and 5.2 from \cite{So4}.

\begin{ex}
(a) Let $p,q\in\zve P\setminus P$ be such that $p\neq q$ and $v(p)=v(q)$. Then $p^2\in\zve(P^2)$ but $qp\in\zve(P^{(2)})$, so $(P^2)\gstr\in v(p^2)\cap\cU\setminus v(qp)$. Since $v(p^2)\neq v(pq)$, this means that $v(x)=v(y)$ does not even imply $v(xz)=_\sim v(yz)$ for $x,y,z\in\zve N$.

(b) Ben De Bondt pointed out that the implication $v(x)=v(y)\Rightarrow v(xz)=v(yz)$ fails even more strongly: for every $x,y\in\zve N$ such that $v(x)=v(y)$ we can find $z\in\zve N$ such that $v(xz)\neq v(yz)$. For example, let $A=\{n!:n\in N\}$ and $f(n)=n!$. Then $z=\zve f(x-1)$ is such that $xz=\zve f(x)\in\zve A$ but $yz\notin\zve A$. Hence $v(xz)\neq v(yz)$.
\end{ex}

So the multiplication in $\zve N$ does not agree with the monad structure: for $\cF,\cG\in\beta N$ the set $\{xy:x\in\mu(\cF),y\in\mu(\cG)\}$ needs not be the monad of an ultrafilter. However, in the next two lemmas we find some regularity in this respect.

\begin{lm}
Let $V(\zve N)$ be ${\goth c}^+$-saturated. For all $\cF,\cG\in\beta N$ the set $M:=\{xy:x\in\mu(\cF),y\in\mu(\cG)\}$ is a union of monads.
\end{lm}

\dokaz We need to prove that $M$ contains all monads that it intersects. So assume $v(z)=v(x_0y_0)$ for some $x_0\in\mu(\cF)$, $y_0\in\mu(\cG)$. The set
$$\Gamma:=\{(x,y)\in\zve N\times\zve N:xy=z\}$$
is internal by The Internal Definition Principle. Moreover, if $A\in\cF$ and $B\in\cG$, then $\Gamma\cap(\zve A\times\zve B)\neq\emptyset$: since $\zve A\zve B=\zve(AB)\in v(x_0y_0)=v(z)$, there are $x\in\zve A$, $y\in\zve B$ such that $z=xy$. It follows that the family $\{\Gamma\}\cup\{\zve A\times\zve B:A\in\cF,B\in\cG\}$ has the finite intersection property. From ${\goth c}^+$-saturation it follows that there are $x\in\mu(\cF)$ and $y\in\mu(\cG)$ such that $z=xy$, so $z\in M$.\kraj

\begin{lm}
Let $\cF,\cG\in\beta N$. $\cF\widemid\cG$ implies $\cF\cdot\cH\widemid\cG\cdot\cH$ and $\cH\cdot\cF\widemid\cH\cdot\cG$ for every $\cH\in\beta N$.
\end{lm}

\dokaz Assume $\cF\cap\cU\subseteq\cG$. Note that, for every $A\in\cU$, every $n\in N$ and every $\cH\in\beta N$, we have $A\subseteq A/n$ and the sets $A/n$ and $A_\cH:=\{m\in N:A/m\in\cH\}$ are in $\cU$ as well.

Let $A\in(\cF\cdot\cH)\cap\cU$, so $A_\cH\in\cF\cap\cU$. It follows that $A_\cH\in\cG$, i.e.\ $A\in\cG\cdot\cH$.

Now let $A\in(\cH\cdot\cF)\cap\cU$, so $A_\cF\in\cH$. Since $A/m\in\cU$ for all $m\in N$, we have $A_\cF\subseteq A_\cG$. Hence $A_\cG\in\cH$, so $A\in\cH\cdot\cG$.\kraj

\section{More on chains}

As we already noted, the $\widemid$-hierarchy gets much more complicated in the "upper" part (above the first $\omega$-many levels). It seems that representing ultrafilters from this area as limits of chains of ultrafilters may be the right way to examine them. In \cite{So4} we considered only chains of length $\omega$, but the following result is a direct generalization to chains of arbitrary order type.

\begin{lm}\label{sublimit}
Every chain $\langle\cF_i:i\in I\rangle$ in $(\beta N,\widemid)$ has the smallest upper bound $\cG_U$ and the greatest lower bound $\cG_L$. Moreover, $\bigcup_{i\in I}(\cF_i\cap\cU)=\cG_U\cap\cU$ and $\bigcap_{i\in I}(\cF_i\cap\cU)=\cG_L\cap\cU$.
\end{lm}

\dokaz We prove the result for the smallest upper bound and the part about the greatest lower bound is proved analogously. If $I$ has the greatest element, this is clearly the wanted bound. Otherwise, let $\cW$ be a nonprincipal ultrafilter on $I$ containing all sets $I_i=\{j\in I:j\geq i\}$ and $\cG_U:=\lim_{i\str\cW}\cF_i$. Let $i\in I$ and $A\in\cF_i\cap\cU$. Since $\langle \cF_i:i\in I\rangle$ is a chain, for every $j\geq i$ we have $A\in\cF_j$ as well. So $\{j\in I:A\in\cF_j\}\supseteq I_i\in\cW$ and thus $A\in\cG_U$. This means that $\bigcup_{i\in I}(\cF_i\cap\cU)\subseteq\cG_U\cap\cU$ and so $\cG_U$ is an upper bound for $\langle\cF_i:i\in I\rangle$.

If we prove $\bigcup_{i\in I}(\cF_i\cap\cU)=\cG_U\cap\cU$ it will imply that $\cG_U$ is the smallest upper bound. So assume there is $A\in\cG_U\cap\cU\setminus\cF_i$ for all $i\in I$. Then $\{j\in I:A^c\in\cF_j\}=I\in\cW$ so $A^c\in\cG_U$, a contradiction.\kraj

If $\cG_U$ is obtained as in the previous proof, we write $[\cG_U]=\lim_{i\in I}\cF_i$. Of course, the bound in question is actually an $=_\sim$-equivalence class, and it follows from the proof that the choice of $\cW$, and thus of the actual representative of the class, is irrelevant.

It is natural to ask: how long can chains in $(\beta N,\widemid)$ be? We answer this for well-ordered chains. (By Example 4.2 from \cite{So4} not all chains are well-ordered.) It is clear that the cardinality of such a chain can not exceed $\goth c$, since there are at most $\goth c$ elements in $\cU$, and in every difference of successive sets there must be one.

In the construction of a chain of maximal cardinality we will use an {\it almost disjoint family}: a family of infinite sets such that every two sets have finite intersection.

\begin{te}
For any $\delta<{\goth c}^+$ there is a chain $\langle\cF_\alpha:\alpha<\delta\rangle$ in $(\beta N,\widemid)$.
\end{te}

\dokaz It is well-known that on every countable set there is an almost disjoint family of cardinality $\goth c$. So let $\{A_\alpha:\alpha<\delta\}$ be such a family of subsets of $P$. For each $\alpha<\delta$ we pick a (prime) ultrafilter $\cP_\alpha\in\overline{A_\alpha}$. Now we define a chain $\langle\cF_\alpha:\alpha<\delta\rangle$ by recursion so that:
\begin{equation}\label{eqniz}
\mbox{each }\cF_\alpha\mbox{ is divisible by all }\cP_\beta\mbox{ for }\beta<\alpha,\mbox{ but does not contain }A_\beta\gstr\mbox{ for }\beta\geq\alpha.
\end{equation}
First let $\cF_0:=1$. Assume $\langle\cF_\beta:\beta<\alpha\rangle$ has been constructed.

Let $\alpha=\beta+1$. We define $\cF_\alpha:=\cF_\beta\cdot\cP_\beta$. Then $\cF_\beta\widemid\cF_\alpha$ (so $\cP_\gamma\widemid\cF_\alpha$ for $\gamma<\beta$ by the induction hypothesis) and $\cP_\beta\widemid\cF_\alpha$. Since, for $\beta\geq\alpha$, neither $\cF_\beta$ nor $\cP_\beta$ contain $A_\beta\gstr$, Lemma \ref{prostilema} implies that $A_\beta\gstr\notin\cF_\alpha$, so $A_\beta\notin\cF_\alpha$. Thus (\ref{eqniz}) holds in this case.

Now let $\alpha\in Lim$. Let $[\cF_\alpha]=\lim_{\beta\in\alpha}\cF_\beta$. By Lemma \ref{sublimit} $\cF_\alpha$ does not contain $A_\beta\gstr$ for $\beta\geq\alpha$ and (\ref{eqniz}) holds again.\kraj

As an immediate corollary we get that $|\cU|={\goth c}$.

We say that $\cG\in\beta N$ is an {\it immediate predecessor} of $\cF\in\beta N$ if $\cG\widemid\cF$, $\cG\neq_\sim\cF$ and there is no $\cH$ such that $\cG\widemid\cH\widemid\cF$ and $\cG\neq_\sim\cH\neq_\sim\cF$.

\begin{lm}
Every ultrafilter $\cF\in\beta N\setminus\{1\}$ either has an immediate predecessor or is the smallest upper bound of a chain of ultrafilters.
\end{lm}

\dokaz We construct a chain $\langle\cG_\alpha\rangle$ of ultrafilters below $\cF$. Let $\cG_0=1$. Assume we have already constructed $\langle\cG_\beta:\beta<\alpha\rangle$. First let $\alpha=\gamma+1$. If $\cG_\gamma$ is an immediate predecessor of $\cF$, we are done. Otherwise let $\cG_\alpha$ be such that $\cG_\gamma\widemid\cG_\alpha\widemid\cF$.

Let $\alpha$ be a limit ordinal. If $\cF\in\lim_{\beta\in\alpha}\cG_\beta$, we are done. Otherwise let $[\cG_\alpha]=\lim_{\beta\in\alpha}\cG_\beta$. This construction ends in less than ${\goth c}^+$ steps, so eventually we get either an immediate predecessor or a desired sequence.\kraj

If the order type of a chain is an ordinal, we can establish another connection with $V(\zve N)$.

\begin{lm}
Let $V(\zve N)$ be ${\goth c}^+$-saturated. For any $\widemid$-chain $\langle\cF_\alpha:\alpha<\gamma\rangle$ there is a $\zvez\mid$-chain $\langle x_\alpha:\alpha<\gamma\rangle$ such that $x_\alpha\in\mu(\cF_\alpha)$ for $\alpha<\gamma$.
\end{lm}

\dokaz We construct the desired sequence by recursion on $\alpha<\gamma$. Let $x_0\in\mu(\cF_0)$ be arbitrary. Assume that $\langle x_\beta:\beta<\alpha\rangle$ has been constructed. For $\alpha=\beta+1$, Lemma \ref{zspo}(a) implies that there is $x_\alpha\in\mu(\cF_\alpha)$ such that $x_\beta\zvez\mid x_\alpha$.

Let $\alpha$ be a limit ordinal. For $\beta<\alpha$ define $\Gamma_\beta=\{y\in\zve N:x_\beta\zvez\mid y\}$. By The Internal Definition Principle each of these sets is internal. We show that the family $F:=\{\zve A:A\in\cF_\alpha\}\cup\{\Gamma_\beta:\beta<\alpha\}$ has the finite intersection property. Let $A_1,\dots,A_k\in\cF_\alpha$ and $\beta_1,\dots,\beta_l\in\alpha$. We may assume that $\beta_1\leq\beta_2\leq\dots\leq\beta_l$, so $x_{\beta_1}\zvez\mid x_{\beta_2}\zvez\mid\dots\zvez\mid x_{\beta_l}$. Applying Lemma \ref{zspo}(a) again we get $y\in\Gamma_{\beta_l}=\Gamma_{\beta_1}\cap\dots\cap\Gamma_{\beta_l}$ such that $y\in\mu(\cF_\alpha)$. But then $y\in\zve A_1\cap\dots\cap\zve A_k$ as well. Thus $F$ has the finite intersection property, so by ${\goth c}^+$-saturation there is $x_\alpha\in\bigcap F$, which concludes the recursion.\kraj

\begin{de}
$\cF\in\beta N$ is an L-limit of ultrafilters if $[\cF]=\lim_{n\in\omega}\cG_n$ for some $\widemid$-chain $\langle\cG_n:n<\omega\rangle$ such that $\cG\in\overline{L_n}$.
\end{de}

It would be nice if every element of the "upper" part of the $\widemid$-hierarchy could be represented as a limit of an $\widemid$-chain of length $\omega$. It would be even better if we could construct such a chain so that all of its elements are from the "lower" part. (An example of this is the representation of the maximal class $MAX=\lim_{n\str\infty}n!$, see \cite{So4}.) Clearly, it would be easy to refine such a representation into a representation as an L-limit. Unfortunately, this is not always possible, as we will see in Lemma \ref{noLlimit}.


\section{$N$-free ultrafilters}

We now encounter another class of ultrafilters beside MAX that is $\widemid$-maximal in certain sense.

\begin{de}
An ultrafilter $\cF\in\beta N$ is $N$-free if it is not divisible by any $n\in N\setminus\{1\}$. A set $A\subseteq N$ is $N$-free if it is an element of some $N$-free ultrafilter.
\end{de}






We call a set $A\subseteq N$ is a {\it strong antichain} if each two distinct $x,y\in A$ are mutually prime. We call it "strong" to distinguish from the notion of an antichain as a set of incomparable elements that we used in \cite{So4}.

\begin{lm}\label{Nfree}
The following conditions are equivalent for any $A\subset N$:

(i) $A$ is $N$-free;

(ii) $A\subseteq n_1N\cup n_2N\cup\dots\cup n_kN$ does not hold for any $n_1,n_2,\dots,n_k\in N\setminus\{1\}$;

(iii) all maximal strong antichains in $A$ are infinite;

(iv) $A$ contains an infinite strong antichain;

(v) $A$ contains arbitrarily long finite strong antichains.
\end{lm}

\dokaz (i)$\Str$(ii) $A\subseteq n_1N\cup n_2N\cup\dots\cup n_kN$ would imply that, for any $\cF\in\overline{A}$, at least one of the sets $n_iN$ is in $\cF$. Thus $\cF$ is not $N$-free, so neither is $A$.\\

(ii)$\Str$(i) The family $F:=\{A\}\cup\{(nN)^c:n\in N\setminus\{1\}\}$ has the finite intersection property (otherwise we would have $A\cap(n_1N)^c\cap\dots\cap(n_kN)^c=\emptyset$, i.e.\ $A\subseteq n_1N\cup\dots\cup n_kN$ for some $n_1,\dots,n_k\in N\setminus\{1\}$). Hence there is an $N$-free ultrafilter containing $A$.\\

(ii)$\Str$(iii) Assume the opposite, that there is a finite strong antichain $X\subseteq A$. Then no $a\in A$ is mutually prime with all elements of $X$. If $P_X$ is the set of all prime divisors of elements of $X$, it follows that $A\subseteq\bigcup_{p\in P_X}pN$, a contradiction with (ii).\\

(iii)$\Str$(iv)$\Str$(v) is obvious.\\

(v)$\Str$(ii) If we assume that $A\subseteq n_1N\cup n_2N\cup\dots\cup n_kN$ for some $n_1,n_2,\dots,n_k\in N$, then $A$ could not contain strong antichains of length more than $k$, since every element of $A$ would be divisible by some of the $n_1,n_2,\dots,n_k$.\kraj

We can now characterize $N$-free ultrafilters in several ways.

\begin{te}
The following conditions are equivalent for every $\cF\in\beta N$:

(i) $\cF$ is $N$-free;

(ii) every $A\in\cF$ contains an infinite strong antichain;

(iii) in every $A\in\cF$ all maximal strong antichains are infinite;

(iv) all maximal strong antichains in $\mu(\cF)$ are infinite;

(v) $\mu(\cF)$ contains two mutually prime elements.
\end{te}

\dokaz (i)$\Str$(ii)$\Str$(iii) follows directly from Lemma \ref{Nfree}.\\

(iii)$\Str$(iv) For every $k\in N$ (iii) implies (taking, for any $A\in\cF$, $B:=A$)
\begin{eqnarray*}
&& (\zs A\in\cF)(\po B\in\cF)(\zs b_1,\dots,b_k\in B)(\po a\in A)\\
&& (\{b_1,\dots,b_k\}\mbox{ is a strong antichain}\Str\{b_1,\dots,b_k,a\}\mbox{ is a strong antichain}).
\end{eqnarray*}
By the bridge theorem \ref{bridge1} it follows that every finite strong antichain in $\mu(\cF)$ can be extended, so every maximal strong antichain in $\mu(\cF)$ is infinite.\\

(iv)$\Str$(v) is trivial.\\

(v)$\Str$(i) Assume the opposite, that $\cF$ is not $N$-free. Then there is $n\in N\setminus\{1\}$ such that $nN\in\cF$. This means that each $x\in\mu(\cF)$ belongs to $\zve(nN)=n\zve N$, so it is divisible by $n$. Hence there are no mutually prime elements in $\mu(\cF)$.\kraj

\begin{te}
There is a $\widemid$-maximal class of $N$-free ultrafilters. (We denote this class by NMAX.)
\end{te}

\dokaz Let $\cU_N$ be the family of all $N$-free sets in $\cU$, and let $\cM=\{(nN)^c:n\in N\setminus\{1\}\}$. $\cU_N$ is closed for intersections: assume the opposite, that $A,B\in\cU_N$ but $A\cap B\notin\cU_N$. Clearly $A\cap B\in\cU$. So by Lemma \ref{Nfree} there is a finite maximal strong antichain $X$ in $A\cap B$. By the same lemma, $X$ is not maximal in either of $A$ or $B$, so there are $a\in A$ and $b\in B$ such that $X\cup\{a\}$ and $X\cup\{b\}$ are strong antichains. However, the least common multiplier $m$ of $a$ and $b$ must be in $A\cap B$, and $X\cup\{m\}$ is also a strong antichain, a contradiction.

We show that $\cU_N\cup\cM$ has the finite intersection property. It suffices to show that every $A\in\cU_N$ intersects $(n_1N)^c\cap\dots\cap(n_kN)^c$ for any $n_1,n_2,\dots,n_k\in N\setminus\{1\}$, but this is exactly what we showed in (i)$\Str$(ii) of Lemma \ref{Nfree}.\kraj

\begin{lm}\label{noLlimit}
No ultrafilter $\cF\in NMAX$ can be represented as an L-limit.
\end{lm}

\dokaz Assume the opposite, that $[\cF]=\lim_{n\in\omega}\cF_n$, with $\cF_n\in\overline{L_n}$. Since an element of $\overline{L_n}$ is divisible by at most $n$ prime ultrafilters, the set $X:=\{\cP\in\overline{P}:(\po n\in\omega)\cP\widemid\cF_n\}$ is countable.

It is not hard to construct an almost disjoint family $\{A_n:n<\omega\}$ of subsets of $P$ such that every $\cP\in X$ belongs to some $\overline{A_n}$. It is well-known that such (countable) family can not be maximal, so there is infinite $B\subseteq P$ almost disjoint from all $A_n$. Since $\cF\cap\cU=\bigcup_{n\in\omega}(\cF_n\cap\cU)$ (Lemma \ref{sublimit}), it suffices to prove that $B\gstr\notin\cF_n$ for all $n$, as this will imply that $\cF\notin NMAX$.

Assume the opposite; then, for some $n<\omega$, the family $\{B\}\cup(\cF_n\cap\cV)$ has the finite intersection property (since the intersection of finitely many sets from $\cF_n\cap\cV$ is still in $\cF_n\cap\cV$, every such set intersects $B\gstr$ and, being in $\cV$, also intersects $B$). Hence there must exist a prime $\cP\in\overline{B}$ such that $\cP\widemid\cF_n$, a contradiction with the definition of $B$.\kraj



\end{document}